\documentclass{article}
\usepackage[leqno]{amsmath}
\usepackage{amsthm,amscd}
\usepackage{amsfonts} 
\usepackage{amssymb} 
\newtheorem{theorem}{Theorem}[section]

\newtheorem{example}[theorem]{Example}
\newtheorem{lemma}[theorem]{Lemma}
\newtheorem{proposition}[theorem]{Proposition}
\newtheorem{corollary}[theorem]{Corollary}
\newtheorem{remark}[theorem]{Remark}

\newcommand{\A}{{\mathcal A}}

\newcommand{\bff}{{\mathbf f}}

\newcommand{\bfx}{{\mathbf x}}
\newcommand{\bfg}{{\mathbf g}}
\newcommand{\sB}{{\mathsf B}}

\newcommand{\sD}{{\mathsf D}}

\newcommand{\bbR}{{\mathbb R}}

\newcommand{\tc}{{\tilde{h}}}

\newcommand{\Der}{{\rm Der}}
\newcommand{\adj}{\operatorname{adj}}

\newcommand{\Jf}{J({\mathbf f})}
\newcommand{\JDkx}{J(D^{k}{\mathbf x})}

\begin{document}
\newcommand{\Vanderdet}{\left|
\begin{matrix}{ccccc}
x_{1} & x_{1}^{3} & . & . &x_{1}^{2l-1}\\
 .    &     .     & . & . &     .      \\
 .    &     .     & . & . &     .      \\
x_{l} & x_{l}^{3} & . & . &x_{l}^{2l-1}     
\end{matrix}
\right|}

\newcommand{\Vanderdetred}{\left|
\begin{matrix}{ccccc}
1 & x_{1}^{2} & . & . &x_{1}^{2l-2}\\
 .    &     .     & . & . &     .      \\
 .    &     .     & . & . &     .      \\
1 & x_{l}^{2} & . & . &x_{l}^{2l-2}     
\end{matrix}
\right|}

\title{Multiderivations of Coxeter arrangements}
\author{{\sc Hiroaki Terao}
\footnote{partially supported by the Grant-in-aid for scientific research (No.1144002), the Ministry of Education, Sports, Science and Technology, Japan
}\\
{\small \it Tokyo Metropolitan University, Mathematics Department}\\
{\small \it Minami-Ohsawa, Hachioji, Tokyo 192-0397, Japan}
}
\date{}
\maketitle
\begin{center} 
{\em Dedicated to Professor Tatsuo Suwa for his sixtieth birthday}
\end{center} 

\medskip

\begin{abstract}
Let $V$ be an $\ell$-dimensional Euclidean space.  
Let $G \subset O(V)$ be a finite
irreducible orthogonal
reflection group.  Let ${\cal A}$ be the corresponding Coxeter
arrangement.  Let $S$ be the algebra of polynomial functions on $V.$
For $H \in {\cal A}$ choose $\alpha_H \in V^*$ such that $H = {\rm
ker}(\alpha_H).$ 
For each nonnegative integer $m$, define the derivation module
$\sD^{(m)}({\cal A}) = \{ \theta \in {\rm Der}_S ~|~
\theta(\alpha_H) \in S \alpha^m_H\}$.
The module is known to be a free $S$-module of rank $\ell$ by K. Saito 
(1975) for $m=1$ and L. Solomon-H. Terao (1998) for $m=2$.
The main result of this paper is that this is the case for all $m$.
Moreover we explicitly construct a basis
for $\sD^{(m)} (\cal A)$.
Their degrees are all equal to $mh/2$ (when $m$ is even) or are equal to
$((m-1)h/2) + m_i (1 \leq i \leq \ell)$ (when $m$ is odd).
Here $m_1 \leq \cdots \leq m_{\ell}$ are the exponents of $G$ and $h=
m_{\ell} + 1$ is the Coxeter number.
The construction heavily uses the primitive derivation $D$ which plays
a central role in the theory of flat generators by K. Saito
(or equivalently the Frobenius manifold structure 
for the orbit space of $G$.)  Some new results concerning 
the primitive derivation $D$  are obtained
in the course of proof of the main result.

\end{abstract}

\bigskip
\setcounter{section}{0}
\setcounter{equation}{0}
\section{Introduction}
~\bigskip

Let $V$ be a Euclidean space of dimension $\ell$ over ${\mathbb R}.$
Let $\cal A$ be a central arrangement of hyperplanes in $V$ :
$\cal A$ is a finite collection of one-codimensional subspaces of $V$.
Let $S$ be the algebra of polynomial functions on $V$.
The algebra $S$ is naturally graded by $S=\oplus_{q \geq 0} S_q $ where
$S_q$ is the space of homogeneous polynomials of degree $q$.
Thus $S_1=V^*$ is the dual space of $V$.
Let ${\rm Der}_S$ be the $S$-module of ${\mathbb R}$-derivations of $S$.
We say that $\theta \in {\rm Der}_S$ is homogeneous of degree $q$ if
$\theta (S_1) \subseteq S_q$.
Choose for each hyperplane $H \in \cal A$ a linear form $\alpha_H
\in V^* $ such that $H=\ker (\alpha_H).$
Let
\begin{eqnarray}
\label{101} 
\sD^{(m)}({\mathcal A})=\{\theta \in {\rm Der}_S ~|~ \theta (\alpha_H) \in 
S\alpha_{H}^{m} {\rm{~for~ any~}} H \in {\cal A} \}
\end{eqnarray}
for each nonnegative integer $m$.
Elements of $\sD^{(m)}(\cal A)$ are called $m$-{\bf derivations} which
were introduced by G. Ziegler \cite {zie}.
Then one has a sequence of inclusions
$$
\sD^{(1)}({\mathcal A}) \supset \sD^{(2)}({\mathcal A}) \supset \cdots .
$$
The arrangement ${\mathcal A}$ is called {\bf free} 
\cite[Chapter 4]{ort} when $\sD^{(1)}
(\mathcal A) = \sD(\mathcal A)$ is a 
free $S$-module.

From now on we assume that $\cal A$ is a Coxeter arrangement : $\cal A$
is the set of reflecting hyperplanes of a finite irreducible subgroup
$G$ of $O(V)$ generated by orthogonal reflections.
Then $G$ naturally acts on $V^*$ and $S$.
Let $S^G$ be the algebra of $G$-invariant polynomials.
Then it is classically known \cite[V.5.3, Theorem 3]{bou} that there
exist algebraically independent homogeneous polynomials $f_1, \cdots, f_{\ell}
\in S^G$ with $\deg f_{1} \leq\dots\leq \deg f_{\ell}$,
which 
are called {\bf basic invariants},
 such that $S^G={\mathbb R}[f_1, \cdots, f_{\ell}].$
Write ${\mathbf f}=(f_{1}, \cdots, f_{\ell})$. 
Let $x_1, \cdots, x_{\ell}$ be a basis for $V^*$ and
 $\partial_i$ be partial derivation with respect to $x_i$:
$\partial_i (x_j)=\delta_{ij}$ (Kronecker's delta).
Let $K$ be the quotient field of $S$.
The {\bf primitive derivation} $D \in {\rm Der}_K$ is characterized
 by
\begin{eqnarray*}
Df_i=\begin{cases}
                                        1 & \text{ for } i=\ell,  \\
                                        0 & \text{ otherwise.}
     \end{cases}
\end{eqnarray*}
The derivation $D$ is, up to a constant multiple,
independent of choice of basic invariants 
$f_1, \cdots, f_{\ell}$ \cite[2.2]{sai2}, \cite[(1.6)]{sys}. 
For ${\mathbf g} = (g_1, \cdots, g_{\ell}) \in K^{\ell} $ 
let $J(\mathbf g)$ be the Jacobian matrix 
with respect to ${\bf x}=(x_1, \cdots, x_{\ell})$:
$J{(\mathbf g)}_{ij}=\partial_i
g_j (1 \leq i, j \leq \ell).$
Let 
$
D^{k} {\mathbf x} = 
(D^{k} x_{1}, \cdots , D^{k} x_{\ell}),
$
where $D^{k} = D\circ D\circ\cdots \circ D$ (composition of $k$ times).
Define $\Gamma\in GL_{\ell}(\bbR)$ by
$\Gamma_{ij} =(x_i, x_j)$, where ( , ) stands for the $G$-invariant
positive definite symmetric bilinear form $V^*$.
The main results of this paper are:

\begin{theorem}
\label{maintheorem} 
For any nonnegative integer
$k$, the $\ell\times\ell$-matrix  
$J(D^{k} {\mathbf x})$ is invertible.  For each nonnegative integer $m$ 
define an $\ell \times \ell$-matrix $P_{m}$ by
$$
               P_{m}  = 
               \begin{cases} 
   \Gamma J(D^{k}{\mathbf x})^{-1}\, & (\text{if } m=2k 
   \text{ is even} ),\\
   \Gamma J(D^{k}{\mathbf x})^{-1} J({\mathbf f})
    & (\text{if } m=2k+1 \text{ is odd}).
               \end{cases}
$$
Define $\xi_1,\ldots,\xi_l \in {\rm Der}_S$ by
$\xi_j  =
\sum_{i=1}^l~ p_{ij} \, \partial_i$ for $1 \leq j \leq l,$ where
$p_{ij}$ is the $(i,j)$-entry of $P_{m}$. 
Then $\xi_i \in \sD^{(m)}({\cal A}),$ and
$\sD^{(m)}({\cal A})$ is a free $S$-module with basis $\xi_1,\ldots,\xi_l.$
Each $\xi_{i} $ 
is homogeneous of degree $kh$ (when $m=2k$) or
of degree $kh+m_{i} $ (when $m=2k+1$) for $i=1,\ldots,\ell$,
where $h$ is the Coxeter number and $m_{1}, \cdots, m_{\ell}$ are the
exponents of $G$.
\end{theorem}

Let $Q := \prod_{H\in\A} \alpha_{H} \in S$.

\begin{theorem}
\label{JDkx} 
Let $k$ be any nonnegative integer.
Then $\det J(D^{k} {\mathbf x})$ is a nonzero multiple of
$Q^{-2k} $. 
\end{theorem}

\begin{corollary}
\label{algindep} 
Let $k$ be any nonnegative integer.
Then the $\ell$ rational functions 
$D^{k} x_{1}, D^{k} x_{2}, \cdots , D^{k} x_{\ell}$
are algebraically independent.
\end{corollary}

\begin{remark}
{\rm For $m=1$, $P_{1} = \Gamma J({\mathbf f})$, $\sD^{(1)}(\A)=\sD(\A)$ 
and Theorem~\ref{maintheorem} was first proved by
K. Saito {\cite[Theorem]{sai1}}, {\cite[Theorem2]{ter}}.
In other words, the Coxeter arrangement is a free arrangement.
For $m=2$, 
Theorem~\ref{maintheorem} was proved in
{\cite[Theorem 1.4]{sot}}.
Theorem~\ref{JDkx} and Corollary~\ref{algindep} for $k=1$ 
were also proved in {\cite[Corollary 3.32]{sot}}.
Corollary~\ref{algindep} for $k=1$ 
had been conjectured by L. Solomon in \cite{sol2} in 1964.}
\end{remark}

\begin{remark}
{\rm The definition (\ref{101}) of $\sD^{(m)}(\A)$ is due to
G. Ziegler \cite[Definition 4]{zie} who developed the theory of
multiarrangements. 
}
\end{remark}

\begin{remark}
{\rm The original motivation to study the module
$\sD^{(m)}(\A)     $  came from the study of the
extended Shi arrangements and the extended Catalan arrangements.
Suppose that $G$ is a Weyl group.   Choose a crystallographic
root system in $V^{*} $ and choose the linear forms 
$\alpha_{H} $ so that $\pm \alpha_{H} $ is a root for each $H\in\A$.
Let $n=\left|\A\right|$ and $\alpha_{1} ,\dots,\alpha_{n} \in V^{*} $  
be a system of positive roots.  
Let  $H_{i, p} $ be the affine hyperplane defined by
$\alpha_{i} = p$ for $p\in{\mathbb Z}$ and $1\leq i\leq n$.
Let $k\geq 1$.  
Define the extended Shi arrangement\cite{shi1}\cite{shi2}\cite{ath2}   
\[
\A^{[-k+1,k]} := \{H_{i,p} \mid 1\leq i\leq n, 
-k+1\leq p\leq k\}, 
\]
and the extended Catalan arrangement
\[
\A^{[-k,k]} := \{H_{i,p} \mid 1\leq i\leq n, -k\leq p\leq k\}. 
\]
For these arrangements 
the characteristic polynomials \cite[p.43]{ort} 
are known \cite{pos1}\cite{pos2}  
\cite[\S 7.1-7.2]{ath}\cite{ath2}  to decompose nicely as:
\begin{eqnarray*}
\chi(\A^{[-k+1, k]}, t) &=& (t - kh)^{\ell},\\ 
\chi(\A^{[-k, k]}, t) &=& \prod_{i=1}^{\ell}  (t - kh-m_{i} ), 
\end{eqnarray*}
when $G$ is of type ${\mathsf{A, B, C, D}}$. 
This result is obtained by direct calculation using the classification
of Weyl groups.  
The same formulas are conjectured to be true for the other
Weyl groups.
We are led to study the module $\sD^{(m)}(\A) $
to understand the meaning of these roots $kh,\dots,kh$ 
(for $\A^{[-k+1, k]} $) 
and $kh+m_{1},\dots,kh+m_{\ell}   $
(for $\A^{[-k, k]} $)
for the characteristic polynomials.   
Let $\A$ be either an extended Shi arrangement or an extended 
Catalan arrangement.
A conjecture 
due to Edelman and Reiner \cite[Conjecture 3.3]{edr}
states that the cone
\cite[p.14]{ort} ${\bf c}\A$ 
 is a free arrangement with exponents
$\{1, kh, \ldots, kh\}$ 
(for $\A^{[-k+1, k]} $) 
or
$\{1, kh+m_{1} , \ldots, kh+m_{\ell} \}$ 
(for $\A^{[-k, k]} $)
\cite[Definition 4.15, Definition 4.25]{ort};
the module 
$D({\bf c}\A)$
is a free module over ${\bf R}[x_{0},\ldots, x_{l}]$.
Athanasiadis \cite{ath} verified this conjecture
for type ${\sf A}_{l} $.  
If the conjecture is true, then, by Ziegler's 
criterion,  which will be stated as
Theorem~\ref{ZieglersTheorem},
we may conclude that 
$\sD^{(2k)}(\A) $  is a free module with generators of degrees
$\{kh,kh,\dots, kh\}$ 
(for $\A^{[-k+1, k]} $) 
and that
$\sD^{(2k+1)}(\A) $  is a free module with generators of degrees
$\{kh+m_{1} ,kh+m_{2} ,\dots, kh+m_{\ell} \}$ 
(for $\A^{[-k, k]} $)
which is true by our main result, Theorem
\ref{maintheorem}.  So Theorem \ref{maintheorem}
may be regarded as a piece of evidence supporting
the conjecture.
}
\end{remark}

%

\section{Notation and preliminaries} 
In this Section we   fix more notation
and state preliminary facts which will be used 
later.
A finite irreducible orthogonal reflection
group  $G$ acts on an $\ell$-dimensional Euclidean space $V$.
Then $G$  acts on the dual space $V^{*}     $  contragradiently.
In other words,
\[
\left<
gx, v
\right>
=
\left<
x, g^{-1} v
\right>
\]
for $x\in V^{*} $, $v\in V$ and $g\in G$ when
$\left<~,~\right>$ denotes the natural pairing
of $V$  and $V^{*}$.
The symmetric algebra $S=S(V^{*})$ 
of $V^{*} $ over $\bbR$ can  be 
identified  as the $\bbR$-algebra of all polynomial functions on
$V$.  The contragradient action of $G$ on $V^{*}$ extends to the action
on $S$ by $(gf)(v) = f(g^{-1} v)$ for $f\in S$, $g\in G$
and $v\in V$.
Each element of the graded $S$-module 
$\Der_{S} $ is an $\bbR$-derivation
of $S$, i.e., an
$\bbR$-linear map
$\theta : S \rightarrow S$ satisfying the Leibniz rule:
$
\theta(f_{1} f_{2} )
=
f_{1} \theta(f_{2})
+
f_{2} \theta(f_{1})
$
for $f_{1}, f_{2} \in S$.
The action of $G$  on $\Der_{S} $  is defined by
$(g\theta) (f) = g(\theta (g^{-1} f))$ for $g\in G$,
$\theta\in \Der_{S} $ and $f\in S$.
The group $G$ naturally acts on 
the quotient field $K$ of $S$ by
\[
g(f_{1} /f_{2} )= (g f_{1})/(g f_{2}) 
\]
for $f_{1} \in S$, $f_{2} \in S \setminus \{0\}$ and  $g\in G$. 
When $f_{1} \in S_{p}$ and $f_{2} \in S_{q} \setminus \{0\}$
define $\deg (f_{1} /f_{2} ) = p - q$. 
For a matrix $A$ of any size with entries in $K$, 
$A_{ij} $ is the $(i, j)$-entry of
$A$.  Let $A^{\top} $ denote the transpose of $A$.
Let $g\in G$. Define the matrix $g[A]$ by 
\[
g[A]_{ij} = g(A_{ij}).
\]

                 Choose a basis $x_{1}, \dots , x_{\ell}  $ for $V^{*}$
and $\bbR$-derivations $\partial_{1},\dots,\partial_{\ell} \in\Der_{S}   $  
such that
\[
\partial_{i} (x_{j}) = \delta_{ij}  
\]
for $1\leq i, j\leq \ell$.
Then $\deg \partial_{i} = 0$ and $\Der_{S} $ is a free $S$-module
with basis $\partial_{1},\dots,\partial_{\ell}  $.
The invertible matrix $\rho(g)\in GL_{\ell}(\bbR) $ is defined by
\begin{eqnarray} 
\label{69} 
g x_{j} = \sum_{i=1}^{\ell} \rho(g)_{ij} x_{i},  
\end{eqnarray} 
or $g[\bfx] = \bfx \rho(g)$. 
Then 
\[
g \partial_{j} = \sum_{i=1}^{\ell} \rho(g)^{-1}_{ji} \partial_{i}.  
\]
One can symbolically express this as
\begin{eqnarray}
\label{58} 
(g \partial_{1}, \dots , g \partial_{\ell}) 
=
\label{gpartial} 
(\partial_{1}, \dots , \partial_{\ell}) (\rho(g)^{\top})^{-1}.
\end{eqnarray}
Let $(~,~)$ denote the positive definite bilinear form 
on $V^{*}$ inverse to the form on $V$.
Recall the $\ell \times \ell$-matrix  
$\Gamma = [(x_{i}, x_{j})]$.
For $g\in G$, we have
\begin{eqnarray}
\label{gGamma}
\Gamma &=& \rho(g)^{\top} \Gamma \rho(g).
\end{eqnarray}
      Choose basic invariants $f_{1}, \dots , f_{\ell}  \in S^{G} $
      with $\deg f_{1} \leq \dots \leq \deg f_{\ell}  $.
      Let $m_{i} = \deg{f_{i} } -1 $ for $i = 1, \dots , \ell$.
      Then the integers $m_{1}, \dots , m_{\ell}  $ are the {\bf exponents}
      of $G$.     
Since $G$ is irreducible, one has   
     \cite[V.6.2]{bou}
     \[
     1=m_{1} < m_{2} \leq\dots\leq m_{\ell-1} < m_{\ell}. 
     \]
     Let $h$  denote the Coxeter number: $h=m_{\ell} + 1 $ 
     \cite[V.6.2, Theorem 1]{bou}. 
Then we have the duality of exponents \cite[V.6.2]{bou}:
\begin{eqnarray*} 
m_{i} + m_{j} = h \,\,\, (\text{if~} i + j = \ell + 1).
\end{eqnarray*}
Write $\bff = (f_{1}, \dots , f_{\ell})$.
Then, for $g\in G$,  we have
\begin{eqnarray}
\label{60}
 g[J(\bff)] = \rho(g)^{-1} J(\bff). 
\end{eqnarray}
By \cite[V.5.5, Proposition 6 (ii)]{bou} we have
\begin{eqnarray}
\label{66}
 \det J(\bff) \doteq Q. 
\end{eqnarray}
Here and elsewhere $\doteq$ stands for equality up to
a nonzero constant multiple.

 Let $\Der_{K} $ be  
 the set of $\bbR$-derivations of $K$.
 The primitive derivation $D$  is defined to be the unique derivation
 satisfying $D f_{i} = 0 \,\, (1\leq i\leq \ell-1)$ and
 $D f_{\ell} = 1$.
 Since 
 \begin{eqnarray*}
 D=\sum_{j=1}^{\ell} (D x_{j}) \partial_{j}   
 \end{eqnarray*}
 we have 
 \[
 \delta_{k\ell} = D f_{k} =  \sum_{j=1}^{\ell} (D x_{j}) (\partial_{j} f_{k}). 
 \]
 Thus $D\bfx = (D x_{1},\dots, D x_{\ell})$ is equal to 
 the bottom row of the matrix $\Jf^{-1} $:
 \begin{eqnarray}
 \label{65} 
 \Jf^{-1}_{\ell j} =  D x_{j}   \,\,\,(j = 1,\dots , \ell).  
 \end{eqnarray}
 By (\ref{66}) we know that $ Q (D x_{j}) \in S$ for $1\leq j \leq \ell$.
 For $k\geq 0$, 
an $\bbR$-linear map $D^{k} : K \rightarrow K$ is defined by
 composing $D$ $k$ times.   Agree that $D^{0} $ is the identity map.
 Then  $D^{k} x_{j} $ has no pole outside the reflecting hyperplanes.
 Since 
 $
 \deg (D x_{j}) = - m_{\ell} = - h + 1
 $
 unless $D x_{j} = 0$, one has
\begin{eqnarray}
\label{degJDkx} 
 \deg (\partial_{i}  D^{k}  x_{j}) &=& - k h
\end{eqnarray}
unless $\partial_{i}  D^{k}  x_{j}=0$.
For $D^{k}\bfx = (D^{k}x_{1},\dots, D^{k}x_{\ell})$
and $g\in G$, we obtain 
$g[D^{k}\bfx] = (D^{k}\bfx) \rho(g)$ and
\begin{eqnarray}
\label{59} 
g[J(D^{k}\bfx)]&=&\rho(g)^{-1}  J(D^{k}\bfx) \rho(g).
\end{eqnarray}
  
\begin{remark}{\rm
The primitive derivation $D$ is, up to a nonzero constant
multiple, independent of choice 
of basic invariants $f_{1},\dots, f_{\ell}$ 
\cite[2.2]{sai2}, \cite[(1.6)]{sys}.  
The derivation $D$ is called the primitive vector field 
by K. Saito in \cite{sai4}
and plays a central role in his theory of flat generators 
and primitive integrals \cite{sai2}
\cite{sai4}.  (See Remark \ref{flat} (ii).) 
}
\end{remark} 

Suppose that another basis $x'_{1},\dots, x'_{\ell}  $ 
for $V^{*}$ is connected with $x_{1},\dots, x_{\ell}$ through
an invertible matrix $M \in GL_{\ell}(\bbR) $:  
  \begin{eqnarray}
  \label{xM}  
  x'_{j} = \sum_{i=1}^{\ell} M_{ij} x_{i},  
  \end{eqnarray} 
  or $\bfx' = \bfx M$. 
  The new objects, which are defined using the new basis
   $x'_{1},\dots, x'_{\ell}  $, will be denoted by 
   $\partial'_{i},   \Gamma', J'$ etc.  
   Then one can symbolically express
\begin{eqnarray}
(\partial'_{1},\dots,\partial'_{\ell}  )
=
(\partial_{1},\dots,\partial_{\ell}  ) (M^{\top})^{-1}  
\label{70} 
\end{eqnarray}

We also have
\begin{eqnarray}
\Gamma' &=& M^{\top} \Gamma M  \label{63} \\
J'(\bff) &=& M^{-1} \Jf \label{62} \\
J'(D^{k} \bfx') &=& M^{-1} J(D^{k} \bfx) M. \label{64} 
\end{eqnarray}

\section{Proof of Theorems} 
In this section we prove Theorems~\ref{maintheorem} and
\ref{JDkx}.  First we split the assertions 
into halves:  

\smallskip

$(E)_{k}$: 
{\it The determinant of $\ell\times\ell$-matrix  
$J(D^{k} {\mathbf x})$ is a nonzero multiple of $Q^{-2k}$ and 
define an $\ell \times \ell$-matrix $P_{2k} $ by
$$
               P_{2k}  = \Gamma J(D^{k}{\mathbf x})^{-1}\, . 
$$
Define $\xi_1,\ldots,\xi_l \in {\rm Der}_S$ by
$\xi_j  =
\sum_{i=1}^l~ p_{ij} \, \partial_i$ for $1 \leq j \leq l,$ where
$p_{ij}$ is the $(i,j)$-entry of $P_{2k}$. 
Then $\xi_j \in \sD^{(2k)}({\cal A}),$ and
$\sD^{(2k)}({\cal A})$ is a free $S$-module with basis $\xi_1,\ldots,\xi_l.$
Each $\xi_{j} $ is homogeneous of degree $kh$ ($1\leq j\leq \ell$)
where $h$ is the Coxeter number.}

\smallskip

$(O)_{k}$:
{\it The determinant of $\ell\times\ell$-matrix  
$J(D^{k} {\mathbf x})$ is a nonzero multiple of $Q^{-2k}$ and 
define an $\ell \times \ell$-matrix $P_{2k+1} $ by
$$
               P_{2k+1}  = \Gamma J(D^{k}{\mathbf x})^{-1} J({\mathbf f})\, . 
$$
Define $\theta_1,\ldots,\theta_l \in {\rm Der}_S$ by
$\theta_j  =
\sum_{i=1}^l~ p_{ij} \, \partial_i$ for $1 \leq j \leq l,$ where
$p_{ij}$ is the $(i,j)$-entry of $P_{2k+1}$. 
Then $\theta_j \in \sD^{(2k+1)}({\cal A}),$ and
$\sD^{(2k+1)}({\cal A})$ is a free $S$-module with basis $\theta_1,\ldots,
\theta_l.$
Each $\theta_{j} $ is homogeneous of degree $kh+m_{j} $ ($1\leq j\leq \ell$)
where $h$ is the Coxeter number and $m_{1}, \cdots, m_{\ell}$ are the
exponents of $G$.}

\smallskip

The following theorem is a special case of a theorem due
to G. Ziegler \cite[Theorem 11]{zie}.
It generalizes
Saito's criterion
\cite[p.270]{sai3}, \cite[Theorem 4.19]{ort}.

\begin{theorem}
\label{ZieglersTheorem}
{\rm\bf (Ziegler's criterion)}
Let $\theta_{1},\ldots , \theta_{\ell} \in D^{(m)}(\A)$.
Then they form a basis for $D^{(m)}(\A)$
if and only if
$\det\left[\theta_{j}(x_{i})\right]_{ij} \doteq Q^{m} $.
\end{theorem}


\begin{lemma}
\label{EkOk}
For any nonnegative integer $k$, the assertion $(O)_{k} $ follows from
$(E)_{k} $.   
\end{lemma}

\begin{proof}
One can symbolically express
\[
(\theta_{1}, \dots , \theta_{\ell})
=
(\partial_{1}, \dots , \partial_{\ell}) P_{2k+1}.
\]
Let $\rho(g)\in GL_{\ell}(\bbR)$ satisfy (\ref{69}):
\[
g x_{j} = \sum_{i=1}^{\ell} \rho(g)_{ij}  x_{i}.    
\]
We have
\begin{eqnarray*}
g[P_{2k+1}]
&=&
g[\Gamma \JDkx^{-1} \Jf]
\\
&=&
\rho(g)^{\top} \Gamma \rho(g)
\rho(g)^{-1} 
\JDkx^{-1}\rho(g) \rho(g)^{-1} \Jf
= 
\rho(g)^{\top} P_{2k+1}
\end{eqnarray*} 
by (\ref{gGamma}), (\ref{59}) and (\ref{60}).
This implies that
each $\theta_{j} $ is $G$-invariant because of (\ref{gpartial}):
\[
(g\partial_{1},\dots,g\partial_{\ell})
=
(\partial_{1},\dots,\partial_{\ell}) (\rho(g)^{\top})^{-1} 
\]
Let $H\in \A$.   Since
\[
\theta_{j}
=
\sum_{i=1}^{\ell}
(\partial_{i} f_{j}) \xi_{i}, 
\]
each $\theta_{j} $ belongs to $\sD^{(2k)}(\A)$.
Write
$\theta_{j} (\alpha_{H}) = \alpha_{H}^{2k} p_{j}(\bfx)  $
for some $p_{j}(\bfx) \in S$.
Let $s_{H}$ be the orthogonal reflection about $H$. 
Then
\begin{eqnarray*} 
\alpha_{H}^{2k}  s_{H}[p_{j}(\bfx)]
&=&
(-\alpha_{H})^{2k}  s_{H}[p_{j}(\bfx)]
=
s_{H}[\alpha_{H}^{2k}  p_{j}(\bfx)]
=
s_{H}[\theta_{j} (\alpha_{H})]\\
&=&
\theta_{j}(s_{H} \alpha_{H}) 
=
\theta_{j} (-\alpha_{H})
 =
-\theta_{j} (\alpha_{H})
=
-\alpha_{H}^{2k}  p_{j}(\bfx). 
\end{eqnarray*} 
Therefore one has 
\[
s_{H}[p_{j}(\bfx)]=
-p_{j}(\bfx).
\]
Hence $p_{j}(\bfx) $ is divisible by $\alpha_{H}$
and thus
$
\theta_{j} (\alpha_{H} ) = 
\alpha_{H}^{2k} p_{j} (\bfx)\in S \alpha_{H}^{2k+1}   $.
This implies $\theta_{j} \in \sD^{(2k+1)} (\A)$
because $H\in \A$ was arbitrary.
Recall Ziegler's criterion~\ref{ZieglersTheorem}.
Since $\det P_{2k+1} = \det (P_{2k} \Jf) \doteq Q^{2k+1}$,
we
can conclude that
$\theta_{1} ,\ldots , \theta_{\ell} \in \sD^{2k+1} ({\A})$ form a basis for the
$S$-module $\sD^{2k+1}({\A})$.  
We also have
\[
\deg \theta_{j} = \deg \xi_{i} + \deg (\partial_{i} f_{j} )
=
kh + m_{j}. 
\]

\end{proof}

\begin{lemma}
\label{Bk}
Let $k$ be any nonnegative integer.
Assume that $(E)_{k} $ holds true.
Define an $\ell \times \ell$-matrix 
$B^{(k+1)}$ with entries in $K$  
by
\begin{equation} 
\label{Bkdef} 
B^{(k+1)}:=-J({\bf f})^{\top}\Gamma J (D^{k+1} {\bf x})J(D^k {\bf x})^{-1}
J({\bf f}).
\end{equation} 
Then

(i) $B^{(k+1)}_{ij}  \in S^G$,

(ii) $\deg B^{(k+1)}_{ij}=m_i+m_j-h$ \,\, $(1\leq i, j\leq \ell)$. 
\end{lemma}

\begin{proof}
(i):
Let $g \in G$.
By (\ref{60}), (\ref{59}) and (\ref{gGamma}), one has
\begin{eqnarray*}
&~& g[B^{(k+1)}]\\
 & = & g[-J({\bf f})^{\top} \Gamma J(D^{k+1}{\bf x})
J(D^k {\bf x})^{-1} J({\bf f})] \\
& = & 
-J({\bf f})^{\top}(\rho(g)^{\top})^{-1}   
\rho(g)^{\top}
\Gamma 
\rho(g)
\rho(g)^{-1} J(D^{k+1}{\bf x})
\rho(g) \rho(g)^{-1} J(D^k {\bf x})^{-1}\rho(g)\rho(g)^{-1} J({\bf f}) \\
& = & B^{(k+1)}.
\end{eqnarray*} 
Thus the matrix $B^{(k+1)}$ is $G$-invariant.

\smallskip
Next we will show that $B^{(k+1)}$ is independent of choice of the basis
$x_1, \cdots, x_{\ell}$ for $V^{*} $.
Suppose that a new basis $x'_1, \cdots, x'_{\ell}$ for $V^{*} $ is connected
with the old basis $x_1, \cdots, x_{\ell}$ through an invertible matrix
$M \in GL_{\ell} ({\mathbb R})$ as (\ref{xM}):
$$
x'_j=\sum ^{\ell}_{i=1} M_{ij} x_i.
$$
The new objects, which are defined using the new basis, will be denoted 
by $\Gamma', J'$ etc.
By (\ref{62}), (\ref{63}) and (\ref{64}), one has
\begin{eqnarray*}
B^{(k+1)'} 
&=& -J'({\bf f})^{\top} \Gamma' J'(D^{k+1}{\bf x'})
J'(D^k {\bf x'})^{-1} J'({\bf f}) \\
&=& -J({\bf f})^{\top}(M^{-1})^{\top}  M^{\top} 
\Gamma M M^{-1}
J(D^{k+1}{\bf x})M
M^{-1}J(D^k {\bf x})^{-1}  M M^{-1}J({\bf f})  \\
&=&  B^{(k+1)}
\end{eqnarray*} 
Therefore $B^{(k+1)}$ is independent of choice of $x_1, \cdots, x_{\ell}$.

Let $H \in {\cal A}$.  Then we may 
 choose an orthonormal basis $x_1=\alpha_H,
x_2, \cdots, x_{\ell}$ for $V^*$ without
affecting $B^{(k+1)}$.
Then $\Gamma$ is the identity matrix.
By Lemma~\ref{EkOk}, the assertion $(O)_k$ holds.
Thus each entry of the first row of $J(D^k {\bf x})^{-1} J({\bf f})$
is divisible by $x^{2k+1}_1$.
Thus, outside the first column, each entry of
$\adj (J(D^k {\bf x})^{-1} J({\bf f}))$ is divisible by $x^{2k+1}_1$.
Since det $J(D^k {\bf x})^{-1} J({\bf f}) \doteq Q^{2k+1}$ 
each entry of 
$$(J(D^k {\bf x})^{-1} J({\bf f}))^{-1} = J({\bf f})^{-1}
J(D^k {\bf x}),
$$
outside the first column, has no pole along 
$x_1=0$.
The $({\ell}, j)$-entry of $J({\bf f})^{-1} J(D^k {\bf x})$ is equal to
$$
\sum^{\ell}_{i=1} J({\bf f})^{-1}_{{\ell}i} \partial_i (D^k x_j)
= \sum^{\ell}_{i=1}(Dx_i) \partial_i (D^k x_j)
=D^{k+1} x_j
$$
by (\ref{65}). 
Thus each entry, except the first, of $D^{k+1}\bfx $
has no pole along $x_{1} = 0$.
Note that the first entry of $D^{k+1}\bfx $ has pole of order
at most $2k+1$ along $x_{1} = 0$.
It follows that each entry of $J(D^{k+1}\bfx)$,
outside the first column, has no pole along $x_{1} = 0$.
Each entry of the first column of $J(D^{k+1}\bfx)$
has pole of order at most $2k+2$ along $x_{1} = 0$.
Since each entry of the first row of $J(D^{k}\bfx)^{-1}\Jf$
is divisible by $x_{1}^{2k+1}  $,
each entry of the product         
 $ J(D^{k+1}\bfx) J(D^{k}\bfx)^{-1}\Jf$
has pole of order at most one along $\alpha_{H} = x_{1} = 0$.
So $$
\alpha_{H} B^{(k+1)} =
-\alpha_{H} 
\Jf^{\top} \Gamma  J(D^{k+1}\bfx) J(D^{k}\bfx)^{-1}\Jf 
$$  
has no pole along $H$.
Since $H\in\A$ was arbitrarily chosen one can conclude
that 
$Q B^{(k+1)}_{ij}  \in S$.
For any $g\in G$, we have
\[
g[Q B^{(k+1)}_{ij} ] = \det (g)^{-1}  Q B^{(k+1)}_{ij}. 
\]
Thus $Q B^{(k+1)}_{ij} $
is an anti-invariant \cite[V. 5,5, Proposition 6 (iv)]{bou} 
and hence  
lies in $Q S^{G} $.
Therefore $B^{(k+1)}_{ij}  \in S^{G}$.
This proves (i).

(ii):
By (\ref{degJDkx}), we have 
\begin{eqnarray*}
\deg B^{(k+1)}_{ij} &=& \deg(
(\partial_{1}  f_{i},\dots,\partial_{\ell}  f_{i})
 \Gamma  
J(D^{k+1}\bfx) J(D^{k}\bfx)^{-1}
(\partial_{1}  f_{j},\dots,\partial_{\ell} f_{j})^{\top} 
) \\
&=& m_{i} - (k+1)h + kh + m_{j} 
=m_{i} + m_{j}- h.   
\end{eqnarray*}
\end{proof} 
  
  \begin{lemma}
  \label{JDg}
  For $\bfg = (g_{1},\dots, g_{\ell})\in K^{\ell} $, we have 
  \smallskip
  
  (i) $J(D\bfg) = J(D\bfx) J(\bfg) + D[J(\bfg)]$, 
  \smallskip
  
  (ii) $D[J(\bfg)^{-1}] = - J(\bfg)^{-1} D[J(\bfg)] J(\bfg)^{-1}$, 
  \smallskip
  
  (iii) $D[J(\bff)] = - J(D\bfx)J(\bff)$,
  
  (iv) $D[J(\bfg)^{-1} J(\bff)] = - J(\bfg)^{-1}J(D\bfg)J(\bfg)^{-1}J(\bff)$. 
  \end{lemma}
  
  \begin{proof} 
(i) The $(i, j)$-entry of $J(D\bfg)$ is equal to:
\begin{eqnarray*}
\partial_{i} (D g_{j})
&=&
\partial_{i} (\sum_{p=1}^{\ell} (D x_{p}) (\partial_{p} g_{j}))
=
\sum_{p=1}^{\ell}  
(\partial_{i} D x_{p}) (\partial_{p} g_{j})
+
\sum_{p=1}^{\ell}  
(D x_{p}) (\partial_{i} \partial_{p} g_{j})\\
&=&
\sum_{p=1}^{\ell}  
(\partial_{i} D x_{p}) (\partial_{p} g_{j})
+
D (\partial_{i} g_{j}),
  \end{eqnarray*}
  which is equal to the $(i, j)$-entry of
  $J(D\bfx) J(\bfg) + D[J(\bfg)]$.
  
  (ii) follows from:
  \[
  0 = 
  D[J(\bfg) J(\bfg)^{-1}]
  =
  D[J(\bfg)] J(\bfg)^{-1}
  +
  J(\bfg) D[J(\bfg)^{-1}].
  \]
  
  (iii) Let $\bfg=\bff$ in (i) to 
  get $J(D\bff)=J(D\bfx)J(\bff) + D[J(\bff)]$.  Apply
  \[
  J(D\bff)=J(Df_{1},\dots ,Df_{\ell}) = J(0,\dots,0,1)=0.  
  \]

  (iv) By (i), (ii), and (iii), one has
  \begin{eqnarray*}
  &~& D[J(\bfg)^{-1} J(\bff)] = - J(\bfg)^{-1}D[J(\bfg)]J(\bfg)^{-1}J(\bff)
  +J(\bfg)^{-1}D[J(\bff)]\\
  &=&
  - J(\bfg)^{-1}\left(J(D\bfg)-J(D\bfx)J(\bfg)\right)J(\bfg)^{-1}J(\bff)
  - J(\bfg)^{-1} J(D\bfx)J(\bff)\\
  &=&
  - J(\bfg)^{-1}J(D\bfg)J(\bfg)^{-1}J(\bff).
  \end{eqnarray*}
  \end{proof} 

\begin{lemma}
\label{Bkdiff}
Let $k$ be any positive integer.
Assume that $(E)_{k-1}$ and $(E)_{k}$ both hold true.
Define 
$B^{(k)} $  
and
$B^{(k+1)} $  
as (\ref{Bkdef}).
Then
\[
B^{(k+1)}- B^{(k)} = B^{(1)} + (B^{(1)})^{\top}.  
\]
  In particular, 
$B^{(k+1)}- B^{(k)} $ is a symmetric matrix.
\end{lemma}

\begin{proof} 
Since 
$B^{(k)}_{ij}  \in S^{G}$ 
and 
\[
\deg B^{(k)}_{ij} = m_{i} + m_{j} - h
\leq 
(h-1)+(h-1)-h = h-2
< \deg f_{\ell}
\]
by Lemma~\ref{Bk}, 
each entry of $B^{(k)}$ lies in 
$\bbR[f_{1},\dots, f_{\ell-1}]$.  
Thus
\begin{eqnarray*}
0&=&D[B^{(k)}]
=
D[-\Jf^{\top} \Gamma J(D^{k}\bfx) J(D^{k-1}\bfx)^{-1}\Jf]\\
&=&
-D[\Jf^{\top}] \Gamma J(D^{k}\bfx) J(D^{k-1}\bfx)^{-1}\Jf
-\Jf^{\top} \Gamma D[J(D^{k}\bfx)] J(D^{k-1}\bfx)^{-1}\Jf\\
&~& -\Jf^{\top} \Gamma J(D^{k}\bfx) D[J(D^{k-1}\bfx)^{-1}\Jf]\\
&=&
-D[\Jf^{\top}] \Gamma J(D^{k}\bfx) J(D^{k-1}\bfx)^{-1}\Jf
-\Jf^{\top} \Gamma D[J(D^{k}\bfx)] J(D^{k-1}\bfx)^{-1}\Jf\\
&~& +\Jf^{\top} \Gamma J(D^{k}\bfx) J(D^{k-1}\bfx)^{-1}
J(D^{k}\bfx)
J(D^{k-1}\bfx)^{-1}\Jf
\end{eqnarray*}
by Lemma~\ref{JDg} (iv). 
Multiply $\Jf^{-1} J(D^{k-1}\bfx) J(D^{k}\bfx)^{-1}\Jf$ 
on the right side 
and we get
\[
0
= - D[\Jf^{\top}]
 \Gamma \Jf -
\Jf^{\top}  \Gamma D[J(D^{k}\bfx)] J(D^{k}\bfx)^{-1}\Jf
- B^{(k)}.
\]
Since
\[
D[\Jf^{\top}] = - \Jf^{\top} J(D\bfx)^{\top}
{\rm ~and~}
D[J(D^{k}\bfx)] = J(D^{k+1}\bfx) - J(D\bfx) J(D^{k}\bfx)
\]
by Lemma~\ref{JDg} (iii) and (i), we finally have 
\[
0 = 
- 
(B^{(1)})^{\top} 
+
B^{(k+1)}
-
B^{(1)}
-
B^{(k)}.
\]
\end{proof} 

\begin{lemma}
\label{Bkplus1}
Let $k$ be a nonnegative integer.
Suppose that 
$(E)_{0}$,
$(E)_{1},\dots,
(E)_{k}$ are all true.
Then 

(i) $B^{(k+1)}  =  (k+1) B^{(1)} + k (B^{(1)})^{\top}$, 

(ii) $\det B^{(k+1)}$ is a nonzero constant. 
   
\end{lemma}

\begin{proof} 
Using Lemma~\ref{Bkdiff} repeatedly we have (i). 

Write $B = B^{(1)} $ for simplicity.
The shape of $B$ is known \cite[Lemma 3.9]{sot}
as follows:  

{\it Case 1)}: If $G$ is not of type ${\sf D}_l$ with $l$ even then
\begin{equation}
\label{Bsoutheast1} 
B =
   \left (
   \begin{matrix}
   0      & 0         &  \cdots  & 0            & B_{1l}    \\
   0      & 0         &  \cdots  & B_{2,l-1}    & B_{2l}    \\
   \vdots & \vdots    &          & \vdots       & \vdots     \\
   0      & B_{l-1,2} &          & B_{l-1,l-1}  & B_{l-1,l} \\
   B_{l1} & B_{l2}    &  \cdots  & B_{l,l-1}  & B_{ll}
   \end{matrix}
   \right )
\end{equation}

\noindent where
\begin{equation}
\label{Bsoutheast2}
\begin{array}{ll}
       B_{ij} = 0  & ~{\rm if}~ i+j < l+1 \\
       B_{ij} \in {\bf R}^* & ~{\rm if}~ i+j = l+1 \\
       m_i B_{ij} = m_j B_{ji} & ~{\rm if}~ i+j = l+1  \, .
\end{array}
\end{equation}

{\it Case 2)}: If $G$ is of type ${\sf D}_l$ with $l = 2p$ then the $2 \times
2$ block in rows and columns $p,p+1$ of the matrix (\ref{Bsoutheast1}) -
the center of the matrix - is to be replaced by a $2 \times
2$ block
\begin{equation*}
B_0 = \left (
   \begin{matrix}
    B_{p,p}    & B_{p,p+1}   \\
    B_{p+1,p}  & B_{p+1,p+1} 
   \end{matrix}
   \right )
\end{equation*}
\noindent with constant entries,
where $B_{p,p+1} = B_{p+1,p}$ and $\det B_0 \in {\bf R}^*.$ 
The statement (\ref{Bsoutheast2}) still holds true outside the 
$2 \times 2$ block $B_{0}$. \medskip

Therefore, by (i), we can verify the following statements:

\medskip
{\it Case 1)}: If $G$ is not of type ${\sf D}_l$ with $l$ even then
\begin{equation*}
\label{B0k} 
\begin{array}{ll}
       B^{(k+1)}_{ij} = 0  & ~{\rm if}~ i+j < l+1 \\
       B^{(k+1)}_{ij} = (k+1) B_{ij} + k B_{ji} = 
       \left(k+1+\frac{k m_{i}}{m_{j}}\right) B_{ij} \in \bbR^{*}  
       & ~{\rm if}~ i+j = l+1  \, .
\end{array}
\end{equation*}

{\it Case 2)}: If $G$ is of type ${\sf D}_l$ with $l = 2p$
then
the $2 \times
2$ block $B_{0}^{(k+1)}  $
 in rows and columns $p,p+1$ of the matrix $B^{(k+1)}$
is 

\begin{equation*}
B^{(k+1)}_0 = \left(
   \begin{matrix}
    B^{(k+1)}_{p,p}    & B^{(k+1)}_{p,p+1}   \\
    B^{(k+1)}_{p+1,p}  & B^{(k+1)}_{p+1,p+1} \\
   \end{matrix}
   \right)=
   (2k+1) B_{0} 
\end{equation*}

\noindent with constant entries,
where $\det B^{(k+1)}_0 \in \bbR^{*} $.
 The statement (\ref{Bsoutheast2}) holds true outside $B^{(k+1)}_{0}$. 
\medskip

In both cases we conclude $\det B^{(k+1)} \in \bbR^{*} $. 
\end{proof}

\begin{remark} 
\label{flat} 
{\rm

(i) 
From what we know about the shape of
$B^{(k)} $ in the proof of 
Lemma \ref{Bkplus1}, we know 
that the symmetric matrix
$B^{(k)}
+
(B^{(k)})^{\top}
$  
has also nonzero constant determinant
for each positive integer $k$.
Since
\[
B^{(1)}
+
(B^{(1)})^{\top}
=
\Jf^{\top} \Gamma D[\Jf]
+D[\Jf]^{\top} \Gamma \Jf
=
D[\Jf^{\top} \Gamma \Jf],
\]
its
$(i,j)$-entry is equal to the Yukawa coupling 
\cite[5.3]{sai4} of $df_{i} $ and 
$df_{j} $. Therefore $B^{(1)}$ can be considered
as an asymmetric ``half'' of the Yukawa coupling
which is symmetric.  The matrices $B^{(k)} $ 
are considered to be higher versions of $B^{(1)} $ 
in the sense that we iteratedely use the primitive 
derivation $D$ to define them.

(ii) Basic invariants 
$f_{1},\dots,f_{\ell}$ are said to be flat generators when
the matrix $B^{(1)}
+
(B^{(1)})^{\top}
=
D[\Jf^{\top} \Gamma \Jf]
$ 
is a constant matrix.
Flat generators always exist \cite{sai2}. 
They give a linear
structure on the quotient
variety $V/G$, which is 
a Frobenius manifold structure
\cite[III\@. \S8]{man}.}
\end{remark} 

{\it Proof of Theorems \ref{maintheorem} and \ref{JDkx}. } 
Thanks to Lemma \ref{EkOk} 
it is enough to verify $(E)_{k} $ for $k\geq 0$.
We prove $(E)_{k} $ by induction on $k$.  
The assertion $(E)_{0}$ is trivially true 
because 
$$\sD^{(0)}({\mathcal A})=\{\theta \in {\rm Der}_S ~|~ \theta (\alpha_H) \in 
S\}=
{\rm Der}_{S}   
$$
and $P_{0}=\Gamma$.  In this case $\xi_{j}(x_{i} ) = (x_{j}, x_{i})$ 
and $\deg \xi_{j} = 0$ for $1\leq j\leq \ell$. 
Suppose that 
$(E)_{0}$,
$(E)_{1},\dots,
(E)_{k}$ are all true.
We will show $(E)_{k+1}$. 
By Lemma \ref{Bkplus1}, $\det B^{(k+1)} $ is a nonzero
constant.   Since
\begin{eqnarray*} 
B^{(k+1)} &=& -J({\bf f})^{\top}\Gamma J (D^{k+1} {\bf x})J(D^k {\bf x})^{-1}
J({\bf f}),\\
\det(J(D^{k}\bfx)) &\doteq& Q^{-2k},\\
\det(\Jf)&\doteq& Q, 
\end{eqnarray*}
we obtain
\begin{eqnarray*}
\det(J(D^{k+1}\bfx)) \doteq Q^{-2k-2}.
\end{eqnarray*}
Let $P_{2k+2} = \Gamma J(D^{k+1}\bfx)^{-1}  $
and $(\xi_{1},\dots,\xi_{\ell}  ) =
(\partial_{1},\dots, \partial_{\ell}  ) P_{2k+2} $.
Then 
\[
\deg \xi_{j} = -\deg(\partial_{i} D^{k+1} x_{j} ) = (k+1)h
\]
for $j=1,\dots,\ell$ by (\ref{degJDkx}).
Hence, in order to verify $(E)_{k+1}$, it suffices to show 
that each $\xi_{j} $ belongs to
$\sD^{2k+2}(\A) $ thanks to
Ziegler's criterion~\ref{ZieglersTheorem}. 
   First we will show how 
$\xi_{1} ,\dots, \xi_{\ell} $ are affected by choosing
a different basis for $V$. Suppose that a new basis
$x'_{1},\dots, x'_{\ell}  $ 
   for $V$  is connected to the old basis
$x_{1},\dots, x_{\ell}  $ 
    through an invertible matrix
    $M\in GL_{\ell}(\bbR) $ as (\ref{xM}). 
     Then, by (\ref{70}), (\ref{63}),
    and (\ref{64}), we obtain 
    \begin{eqnarray*}
    (\xi'_{1},\dots,\xi'_{\ell})
    &=&(\partial'_{1},\dots,\partial'_{\ell})\Gamma' J'(D^{k+1}\bfx')^{-1}\\
    &=&(\partial_{1},\dots,\partial_{\ell}) (M^{\top})^{-1}   M^{\top} 
    \Gamma M  M^{-1}  J(D^{k+1}\bfx )^{-1} M  \\
    &=&(\xi_{1},\dots,\xi_{\ell}) M. 
    \end{eqnarray*}
    In other words,
    $\xi_{1},\dots,\xi_{\ell}$ satisfy the same base change rule as
    $\bfx = (x_{1},\dots, x_{\ell})$.
    
    Let $H\in\A$.  Then we may assume that $x_{1} = \alpha_{H}, x_{2} ,\dots,
    x_{\ell}  $ are an orthonormal basis.   It is enough to show that
    $\xi_{j}(x_{1} ) \in S x_{1}^{2k+2}  $ for each $j$.
    Since $(O)_{k} $ follows from $(E)_{k} $ by Lemma~\ref{EkOk},
    each entry of the first row of $J(D^{k+1}\bfx )^{-1}\Jf$ 
    is divisible by $x_{1}^{2k+1}  $.
    By the exactly same argument as in the proof of Lemma~\ref{Bk},
    we know that each entry of $J(D^{k+1}\bfx )$, except the first column,
    has no pole along $x_{1} = 0$.
    Thus each entry of the first row of $\adj J(D^{k+1}\bfx )$ has no pole 
    along $x_{1} =0$.  Since $\det J(D^{k+1}\bfx ) \doteq Q^{-2k-2} $,
    each entry of the first row of $J(D^{k+1}\bfx )^{-1}$ is divisible by
    $x_{1}^{2k+2}  $.  This implies   
          $\xi_{j}(x_{1} ) \in S x_{1}^{2k+2}  $ for each $j$.
          Thus $(E)_{k+1} $ holds true and the induction proceeds.
          This completes the proof of Theorems \ref{maintheorem} and
          \ref{JDkx}.  
    
 \medskip
 
 \begin{example}
 {\rm {\bf (The case $\sB_{2}$)}
 Let us explicitly calculate $J(D\bfx), P_{3}, \xi_{1}^{(3)}, 
 \xi_{2}^{(3)}$ when $\A$ is of type 
 $\sB_{2}$. Let $x_{1}$ and $x_{2}$ be an orthonormal basis and  
 \[
 f_{1} = \frac{1}{2}  (x_{1}^{2} + x_{2}^{2}),
 \,\,\,
 f_{2} = \frac{1}{4} (x_{1}^{4} + x_{2}^{4}).
 \]
 Then
\[
P_{1} = J({\bf f}) =
\left(
\begin{matrix}
x_{1}  & x_{1}^{3}\\
x_{2}  & x_{2}^{3}
\end{matrix}
\right).
\]
 Choose $Q = \det \Jf = x_{1} x_{2} (x_{1} +x_{2} )(x_{2} - x_{1}).$
 Since
 \[
 \Jf^{-1} = 
 \frac{1}{Q}  
 \left(
\begin{matrix}
x_{2}^{3}   & -x_{1}^{3}\\
-x_{2}  & x_{1}
\end{matrix}
\right),
\]
one has
\[
D = -\frac{x_{2}}{Q} \partial_{1} + \frac{x_{1}}{Q} \partial_{2},
\,\,\,
D\bfx = -\left(
\frac{x_{2}}{Q}, \frac{x_{1}}{Q}\right),
\]
\[
J(D\bfx)
=
\frac{1}{Q^{2}}
\left(
\begin{matrix}
(3 x_{1}^{2} - x_{2}^{2})x_{2}^{2}   & 2x_{1}^{3} x_{2}  \\
2 x_{1} x_{2}^{3}   & x_{1}^{2}(x_{1}^{2} - 3 x_{2}^{2}) 
\end{matrix}
\right),
\]
\begin{eqnarray*} 
P_{3} 
&=&
J(D\bfx)^{-1} \Jf\\ 
&=&
\left(
\begin{matrix}
-(1/3) x_{1}^3 (x_{1}^2 - 5 x_{2}^2) 
& -(1/3) x_{1}^3 (x_{1}^4 - 3 x_{1}^2 x_{2}^2 - 2 x_{2}^4)\\
   (1/3) x_{2}^3 (5 x_{1}^2 - x_{2}^2) & 
   (1/3) x_{2}^3 (2 x_{1}^4 + 3 x_{1}^2 x_{2}^2 - x_{2}^4)
\end{matrix} 
\right).
\end{eqnarray*} 
Thus the basis $\xi^{(3)}_{1}, \xi^{(3)}_{2}$ for $D^{(3)}(\A)$ 
are given by 
\begin{eqnarray*} 
\xi^{(3)}_{1} &=& 
-(1/3) x_{1}^3 (x_{1}^2 - 5 x_{2}^2) \partial_{1} 
+ (1/3) x_{2}^3 (5 x_{1}^2 - x_{2}^2) \partial_{2},\\
\xi^{(3)}_{2} &=& 
-(1/3) x_{1}^3 (x_{1}^4 - 3 x_{1}^2 x_{2}^2 - 2 x_{2}^4)\partial_{1} 
+ (1/3) x_{2}^3 (2 x_{1}^4 + 3 x_{1}^2 x_{2}^2 - x_{2}^4)\partial_{2}.
\end{eqnarray*} }
 \end{example} 

 \bigskip
 
 We have the following alternative inductive expression for the matrices
 $P_{m} $ in Theorem~\ref{maintheorem} using the matrices $B^{(k)}$:
 
 \begin{proposition}
 \label{PB} For any nonnegative integer $m$, we have $$
 P_{m} =
 \begin{cases}
 \Gamma     & (\text{if~} m = 0),\\
 P_{m-1} \Jf & (\text{if~} m \text{~is 
 odd}),\\
 - P_{m-1} (B^{(m/2)})^{-1} P_{1}^{\top} & (\text{if~} m \text{~is 
 positive and even }).
  \end{cases}
  $$
 \end{proposition}
  
 \begin{proof}
 Only the third formula needs to be verified.
 Let $m=2k$. Compute
 \begin{eqnarray*} 
 &~&-P_{2k-1} (B^{(k)})^{-1} P_{1}^{\top}\\
 &=& - \Gamma J(D^{k-1}\bfx)^{-1} \Jf (-\Jf^{\top} \Gamma \JDkx  
 J(D^{k-1}\bfx)^{-1} \Jf)^{-1}\Jf^{\top}\Gamma\\
 &=& \Gamma \JDkx^{-1} = P_{2k}.  
 \end{eqnarray*} 
\end{proof} 

\begin{example}{\rm {\bf (The case $\sB_{\ell}$)}
Let $G$ be the Coxeter group of type ${\sf B}_{\ell}$
acting on an $\ell$-dimensional
Euclidean
space $V$
by signed permutations of an orthonormal basis
$e_{1}, \ldots, e_{\ell}.$
Let $x_{1}, \ldots, x_{\ell}$ be the dual basis
for $V^{*} $.  Then $\Gamma$ is the identity matrix.
Define
$$
f_{i} = f_{i}(x_{1}, \ldots, x_{\ell})
=
\frac{1}{2i} \sum_{j=1}^{\ell} x_{j}^{2i}      
$$
for $1\leq i \leq \ell$.   
We will use the basic invariants
$f_{1},\dots,f_{\ell}$. 
Then
$$
P_{1} =
J({\bf f}) =
\left(
\begin{matrix}
x_{1}  & x_{1}^{3} & \cdots &x_{1}^{2\ell-1}\\
\vdots & \vdots    &        & \vdots     \\
x_{\ell}  & x_{\ell}^{3} & \cdots &x_{\ell}^{2\ell-1}
\end{matrix}
\right).
$$
For $i \geq 0$ let
$$
h_{i} (x_{1},\ldots, x_{\ell} )
=
\sum_{i_{1} + \cdots + i_{\ell} = i} x_{1}^{i_{1} } \cdots x_{\ell}^{i_{\ell} }
$$
be the $i$-th complete symmetric polynomial.
Let
$$
\tc_{i}
=
 \tc_{i} (x_{1},\ldots, x_{\ell} )
=
h_{i} (x_{1}^{2} ,\ldots, x_{\ell}^{2}  ) \, .
$$
Then $\tc_{i} $ is a $G$-invariant polynomial of degree
$2i$.   Define
$\tc_{i} (x_{1},\ldots, x_{\ell} ) = 0$ if
$i < 0$.
In \cite[\S 5.2]{sot}, $B^{(1)} $ was 
calculated 
as:
\[
B^{(1)}_{ij} = 
(2j-1) \tc_{i+j-\ell-1}(x_{1},\ldots, x_{\ell})     
\]
for $1\leq i,j\leq\ell$.
Then, by Lemma \ref{Bkplus1},
we have
\[
B^{(k)}_{ij} = k B^{(1)}_{ij} + (k-1) B^{(1)}_{ji}
=
\left\{k(2i+2j-2)-2i+1\right\} \tc_{i+j-\ell-1}(x_{1},\ldots, x_{\ell})     
\]
for $1\leq i,j\leq\ell$ and $
k\geq 1$.  We can inductively calculate $P_{m} $ using Proposition
\ref{PB}.  
   
}
\end{example}

{\small {\bf Acknowledgment}: The author would like to
express his gratitude to the referee
for carefully reading the manuscript and
giving him valuable suggestions.}

\end{document}